\documentclass[11pt]{article}
\usepackage{amssymb,amsmath,latexsym}
\usepackage{epsfig}

\allowdisplaybreaks

\newtheorem{thm}{Theorem}[section]

\newtheorem{lemma}[thm]{Lemma}

\newtheorem{defn}[thm]{Definition}

\newtheorem{remark}[thm]{Remark}

\numberwithin{equation}{section}

\def\pf{{\medskip\noindent {\bf Proof. }}}
\def\qed{{\hfill $\Box$ \bigskip}}

 \def\sB {{\cal B}} 
 \def\sE {{\cal E}} \def\sF {{\cal F}}
 \def\sH {{\cal H}} 
  \def\sL {{\cal L}}

 \def\bE {{\mathbb E}}

\def\bP {{\mathbb P}}  \def\bR {{\mathbb R}}

\def\1{{\bf 1}}

\def\wt{\widetilde}

\def\E{{\mathbb E}}
\def\P{{\mathbb P}}

\def\<{\langle}
\def\>{\rangle}
\def\Exp{{\rm Exp}}
\def\loc{{\rm loc}}

\def\K{{\bf K}}
\def\J{{\bf J}}

\def\bea{\begin{align*}}
\def\eea{\end{align*}}
\def\bee{\begin{equation}}
\def\eee{\end{equation}}

\def\eps{\varepsilon}

\makeatletter
\@addtoreset{equation}{section}

\makeatother

\begin{document}
\bibliographystyle{plain}

\title{\Large \bf
$L^p$-independence of
spectral bounds of generalized non-local Feynman-Kac semigroups}

\author{{\bf Zhen-Qing Chen}\thanks{Research partially supported
by NSF Grants DMS-0906743 and DMR-1035196.}}
\date{(February 8, 2012)}
\maketitle

 \begin{abstract}
Let $X$ be a symmetric strong Markov process on a Luzin space.
In this paper,
we present criteria of the $L^p$-independence of spectral bounds for
generalized non-local Feynman-Kac semigroups of
$X$ that involve continuous additive functionals of $X$ having
zero quadratic variations and discontinuous additive functionals of $X$.
 \end{abstract}



 \section{Introduction}

Transformation by multiplicative functionals is one of the most important transforms for Markov processes. Feynman-Kac transforms
and Girsanov transforms are particular cases.
They  play an important role in the probabilistic
as well as analytic aspect of   potential theory.
See \cite{CFKZ1, CFKZ2, CZ2} and the references therein
for some of the recent results in the context of symmetric Markov processes.

Suppose that  $E$ is  a Lusin space
(i.e., a space that is homeomorphic to a Borel subset of a
compact metric space)
and $\sB(E)$ denotes the Borel $\sigma$-algebra on $E$.
Let $m$ be a Borel $\sigma$-finite measure on $E$
with $\hbox{supp} [ m ] =E$ and
$X =(\Omega,\, \sF,\, \sF_t, \, X_t, \, \P_x, \, x\in E)$ be an
$m$-symmetric irreducible
  Borel standard process on $E$ with lifetime $\zeta$
(cf. Sharpe \cite{sharpe} for the terminology). For a continuous additive
functional $A$ of $X$ having finite variation,
one can do Feynman-Kac transform:
$$ T_t f (x)= \E_x \left[ e^{A_t} f(X_t)\right], \qquad t\geq 0.
$$
It is easy to check (see \cite{ABM}) that, under suitable Kato class condition on $A$,
$\{T_t; t\geq 0\}$ forms a strongly continuous symmetric semigroup on $L^p(E; m)$
for every $1\leq p\leq \infty$ and that its $L^2$-infinitesimal generator
is $\sL^\mu:=\sL + \mu$, where $\sL$ is the $L^2$-infinitesimal generator  of the process $X$
and $\mu$ is the (signed) Revuz measure for the continuous additive functional
$A$. To emphasize the correspondence between continuous additive functionals
and Revuz measures, let's denote $A$ by $A^\mu$.
In fact, the process $X$ has many continuous additive functionals that do not have finite variations. For example, for $u$ in the extended Dirichlet space $\sF_e$ of $(\sE, \sF)$, $u(X_t)-u(X_0)$ has
 the Fukushima's decomposition $M^u+N^u$, where $M^u$ is a square integrable martingale additive functional of $X$ and $N^u$ is a continuous additive
function that in general is only of zero quadratic variation.
It is also natural to consider the following
 generalized Feynman-Kac transform:
$$ T_t f(x)= \E_x \left[ e^{N^u_t} f(X_t)\right], \qquad t\geq 0.
$$
Let $\mu_{\< u\>}$ be the Revuz measure for $\<M^u\>$, the quadratic
variation process of $M^u$. We refer the reader to the Introduction
of \cite{CZ} for a brief history of the above transformation by $N^u$.
It is shown in \cite{CZ} that when $\mu_{\< u\>}$ is in Kato class of $X$,
$\{T_t; t\geq 0\}$ forms a strongly continuous symmetric semigroup on
$L^2(E; m)$ and its associate quadratic form is $(Q, {\cal D}(Q))$, where
${\cal D}(Q)_b\subset \sF_b$ and
$$ Q(f, g)= \sE (f, g) + \sE (fg, u) \qquad \hbox{for } f, g\in \sF_b.
$$
Here for a function space $\sH$, we use $\sH_b$ to denote space of
bounded functions in $\sH$.
When the process $X$ is discontinuous, it has many discontinuous additive functionals.
Let $F$ be a bounded symmetric function on $E\times E$
that vanishes along the diagonal $d$
of $E\times E$. We always extend it to be zero off $E\times E$.
Then $\sum_{0<s\leq t} F(X_{s-}, X_s)$, whenever it is summable,
 is an additive functional of $X$. Hence one can perform generalized
 non-local Feynman-Kac transform
 \begin{equation}\label{e:1.1}
  T^{u, \mu, F}_t f (x):= \E_x \left[ \exp \left(N^u_t+ A^\mu_t+\sum_{0<s\leq t} F(X_{s-}, X_s) \right) f(X_t)\right],
   \qquad t\geq 0.
\end{equation}
Under some suitable Kato class conditions on the  measures $\mu_{\<u\>}$, $\mu$ and the function $F$, it can be shown (see Theorem \ref{T:3.5} below) that
$\{T^{u, \mu, F}_t; t\geq 0\}$ is a strongly continuous symmetric semigroup on
$L^p(E; m)$ for every $1\leq p\leq \infty$. Hence the limit
$$ \lambda_p (X; u+ \mu+F):= - \lim_{t\to \infty} \frac1t \log \| T^{u, \mu, F}_t\|_{p, p}
$$
exists, which will be called the $L^p$-spectral bound of the
generalized non-local Feynman-Kac
semigroup $\{T^{u, \mu, F}_t; t\geq 0\}$.
We will show in this paper that under suitable  conditions,
 $\lambda_p (X; u+\mu+F) = \lambda_2 (X; u+\mu+F)$
for all $1\leq p\leq \infty$ if $\lambda_2 (X; u+\mu+F)\leq 0$.
If in addition $X$ is conservative, then $\lambda_2 (X; u+\mu+F)\leq 0$
becomes a necessary and sufficient condition for the independence of
$\lambda_p (X, u+\mu+F)$ in $p\in [1, \infty]$.
The $L^2$-spectral bound $\lambda_2 (X; u+\mu+F)$ has a
 variational formula in terms of the Dirichlet form of $X$, $\mu$ and $F$, see
\eqref{e:3.9} below.

When $F=0$ and $u=0$, the $L^p$-independence of spectral bounds
for  continuous Feynman-Kac transforms $\{T^{0, \mu, 0}_t, t\geq 0\}$
was investigated by
Takeda in \cite{T, T2} for conservative Feller processes or symmetric Hunt processes satisfying strong Feller property and a
tightness assumption, respectively,
both using a large deviation approach.
The results in \cite{T} were extended
to purely discontinuous Feynman-Kac transforms
$\{T^{0, 0, F}_t, t\geq 0\}$
(i.e. with $u=0$ and $\mu=0$) first
in \cite{TT} for rotationally symmetric $\alpha$-stable processes
and then in \cite{Taw} for   conservative
doubly Feller processes,
both papers  again using a large deviation approach similar
to those in \cite{T, T2}.
A stochastic process is said to be doubly Feller
if it is a Feller process that has the strong Feller property.
The $L^p$-independence of spectral bounds
for continuous generalized Feynman-Kac transforms $\{T^{u, \mu, 0}_t, t\geq 0\}$ (i.e. with $F=0$)
is studied recently in \cite{DKK} for doubly Feller processes
on a locally compact metric space $E$
and for those $u\in \sF_e$ that is continuous on $E$ and vanishes at infinity,
 also using a large deviation approach refined
from \cite{T, T2}.
 In a very recent paper \cite{C3} by the author,   a completely
  different approach is developed to study the $L^p$-independence
 of spectral bounds for non-local Feynman-Kac semigroups $\{T^{0, \mu, F}_t, t\geq 0\}$ (i.e. with $u=0$)
for symmetric Markov processes that may not have strong Feller property,
using the gaugeability results established in \cite{C} for
continuous Feynman-Kac functionals.
This new approach yields new criteria
for the $L^p$-independence of spectral bound even for local Feynman-Kac
semigroups.

The approach of this paper is  different from that of \cite{DKK}.
We do not use large deviation theory.
Using the idea from \cite{CZ},  we decompose
transformation by multiplicative
functional $e^{N^u_t}$ into a combination of a Girsanov transform,
a continuous Feynman-Kac transform  followed by an $h$-transform.
So essentially, after a Girsanov transform, we can reduce
the generalized non-local Feynman-Kac transform into
a non-local Feynman-Kac transform for a new process.
We can  then apply the criteria from \cite{C3} to the latter
to obtain criteria of the $L^p$-independence of spectral radius
for generalized Feynman-Kac semigroup $\{T^{u, \mu, F}_t, t\geq 0\}$.

To keep the exposition of this paper as transparent as possible,
 we have not attempted to present the most general conditions
 on $u$, $\mu$ and $F$. For example, by applying results from
 \cite{C3} for continuous Feynman-Kac transforms instead of that for
 non-local Feynman-Kac transforms, conditions on $\mu$ can
 be weakened for $\{T^{u, \mu, 0}_t, t\geq 0\}$
 in the case of  $F=0$.

  The rest of the paper is organized as follows.
  In Section \ref{S:2}, we give precise setup of this paper, including the
  definitions of Kato classes and L\'evy systems and recalling
  the main results from \cite{C3} that will be used in the sequel.
    Generalized non-local Feynman-Kac transform and its reduction to non-local Feynman-Kac transform via Girsanov transform are studied in
    Section \ref{S:3}. The criteria of the $L^p$-independence of
     spectral bound for generalized non-local Feynman-Kac semigroups
     are established in Section \ref{S:4}. Several examples are given
     in Section \ref{S:5} to illustrate the main results of this paper.

\section{Kato classes and non-local Feynman-Kac transform}\label{S:2}

Let $E$ be a Lusin space
and $\sB(E)$ be the Borel $\sigma$-algebra on $E$.
Let $m$ be a Borel $\sigma$-finite measure on $E$
with $\hbox{supp} [ m ] =E$ and
$X =(\Omega,\, \sF,\, \sF_t, \, X_t, \, \P_x, \, x\in E)$ be an
$m$-symmetric irreducible
transient Borel standard process on $E$ with lifetime $\zeta$.
We like to point here that, since we are only concerned with
 the Schr\"odinger semigroups of $X$,
the transience assumption on $X$
is just a matter of convenience and is unimportant---we can always consider the 1-subprocess $X^{(1)}$ of $X$ instead of $X$ if necessary.
Let $(\sE, \sF)$ denote the Dirichlet form of $X$;
that is, if we use $\sL$ to denote the infinitesimal generator of $X$,
then $\sF$ is the domain of the operator $\sqrt{-\sL}$ and
for $u, v \in \sF$,
$$ \sE (u, v) = (\sqrt{-\sL u}, \, \sqrt{-\sL v} )_{L^2(E; m)}.
$$
We refer the reader to \cite{CF} or \cite{FOT} for terminology and various
properties of Dirichlet forms such as continuous additive functional,
martingale additive functional, extended Dirichlet space.

The transition operators $P_t$, $t\ge 0$, are defined by
$$
P_tf(x):=\E_x [f(X_t)] =\E_x[ f(X_t); \, t<\zeta ].
$$
(Here and in the sequel, unless mentioned otherwise,
we use the
convention that a function defined on $E$ takes the value 0 at the
cemetery point $\partial$.)
We assume that there is a Borel symmetric function $G(x, y)$
on $ E\times E$ such that
$$
\E_x \left[ \int_0^\infty f(X_s) ds \right]
=\int_E G(x, y) f(y) m(dy)
$$
for all measurable $f\geq 0$.
$G(x, y)$ is called the Green function of
$X$. The Green function $G$ will always be chosen so that for each fixed
$y\in E$,  $x\mapsto G(x, y)$ is an
excessive function of $X$.
This choice of the Green function is always possible; see \cite{sharpe}.

For every $\alpha>0$, one deduces from
the existence of the Green function $G(x, y)$
that there exists
a kernel $G_\alpha(x, y)$ so that
$$
\E_x \left[ \int_0^\infty e^{-\alpha s} f(X_s) ds \right]
=\int_E G_\alpha (x, y) f(y) m(dy), \qquad x\in E,
$$
for all measurable $f\geq 0$. Clearly, $G_\alpha (x, y)\leq G(x, y)$.
Note that by \cite[Theorem 4.2.4]{FOT}, for every $x\in E$ and $t>0$,
$X_t$ under $\P_x$ has a density function $p(t, x, y)$
with respect to the measure $m$.

A set $B$ is said to be $m$-polar if
$\P_m (\sigma_B<\infty )=0$, where
$\sigma_B:=\inf \{ t>0: \, X_t \in B \}$.
We call a positive measure $\mu$ on $E$ a smooth measure
of $X$
if there is a positive continuous additive functional
(PCAF in abbreviation)
$A$ of $X$ such that
\begin{equation}\label{eqn:Revuz1}
\int_E f(x) \mu (dx) = \uparrow \lim_{t\downarrow 0}
\E_m \left[ \frac1t \int_0^t f(X_s) dA_s \right].
\end{equation}
for any Borel $f\geq 0$.
Here $\uparrow \lim_{t\downarrow 0}$ means the quantity
is increasing
as $t\downarrow 0$.
The measure $\mu$ is called the
Revuz measure of $A$.
We refer to \cite{CF, FOT}
 for the characterization of smooth
measures in terms of nests and capacity.

For any given positive smooth measure $\mu$, define
$G \mu (x) =\int_E G(x, y) \mu (dy)$.
It is known (see Stollmann and Voigt \cite{SV})
that for any positive smooth measure
$\mu$ of $X$,
\begin{equation}\label{eqn:en}
\int_E u(x)^2 \mu(dx) \leq \|G\mu \|_\infty \,
\sE (u, u) \qquad \hbox{for } u \in \sF.
\end{equation}
Recall that as $X$ is assumed to have a Green function,
any $m$-polar set is polar. Hence by
(\ref{eqn:en}) a PCAF $A$ in the sense of
\cite{FOT} with an exceptional set that has a bounded potential
(that is, $x \mapsto \E_x \left[ A_\zeta \right]=G \mu$ is bounded
almost everywhere on $E$,
where $\mu$ is the Revuz measure of $A$)
can be uniquely
refined into a PCAF in the strict sense (as defined on p.195 of
\cite{FOT}).
This can be proved by using the same argument as that in the
 proof of Theorem 5.1.6 of   \cite{FOT}.

\medskip

For a signed measure $\mu$, we use $\mu^+$ and $\mu^-$ to denote
the positive part and negative part of $\mu$  appearing in
the Hahn-Jordan decomposition of $\mu$.
The following definitions are taken from Chen \cite{C}.

\begin{defn}\label{D:2.1}
Suppose that $\mu$ is a signed smooth measure.
Let $A^\mu$ and $A^{|\mu|}$ be the continuous additive
functional and positive continuous additive functional
 of $X$ with Revuz measures
 $\mu$ and $|\mu |$, respectively.

\begin{description}
\item{\rm (i)} We say  $\mu$ is
in the Kato class of $X$, $\K(X)$ in abbreviation,
if
$$ \lim_{t\to 0} \sup_{x\in E}
\E_x \left[ A^{|\mu |}_t \right] =0.
$$

\item{\rm (ii)} $\mu$  is said to be in the class $\K_\infty (X)$
if for any $\eps>0$, there is a Borel set
$K=K(\eps)$ of finite $|\mu|$-measure
and a constant $\delta=\delta(\eps)>0$
such that for all measurable set $B\subset K$
with $|\mu|(B)<\delta$,
\begin{equation}\label{eqn:1.1}
\| G (\1_{K^c\cup B} | \mu |) \|_\infty <\eps.
\end{equation}

\item{\rm (iii)} $\mu$ is said to be in the class $\K_1 (X)$
if there is a Borel set $K$ of finite $|\mu|$-measure
and a constant $\delta>0$
such that
\begin{equation}\label{eqn:1.3}
\beta_1 (\mu):= \sup_{B\subset K: \, |\mu|(B)<\delta}
\| G (1_{K^c \cup B} | \mu |) \|_\infty <1.
\end{equation}

\item{\rm (iv)} A function $q$ is said to be in class
$\K (X)$,  $\K_\infty (X)$ or $\K_1 (X)$
 if $\mu(dx) :=q(x) m(dx)$ is in the corresponding spaces.
\end{description}
\end{defn}

According to \cite[Proposition 2.3(i)]{C}, $\K_\infty (X)\subset \K (X)\cap \K_1(X)$.
Suppose that $\mu$ is a positive measure in $\K_1 (X)$.
By Propositions 2.2 in \cite{C}, $G\mu (x)=\E_x [ A^\mu_\infty]$
 is bounded
and so (\ref{eqn:en}) is satisfied. Therefore the PCAF
corresponding to $\mu$ can and is always taken to be
in the strict sense.

\medskip

Let $(N, H)$ be a L\'evy system for $X$
(cf. Benveniste and Jacod \cite{BJ}
and Theorem 47.10 of Sharpe \cite{sharpe});
that is,
$N(x, dy)$  is a kernel from $(E , {\cal B}(E))$ to
$(E, {\cal B}(E))$ satisfying
$N(x, \{ x \})=0$,
and $H_t$ is a PCAF of $X$ with bounded
1-potential such that  for any nonnegative Borel function $f$ on $E \times
E$ vanishing on the diagonal and any $x\in E$,
\begin{equation}\label{eqn:LS1}
\E_x\left[ \sum_{s\le t}f(X_{s-}, X_s)\1_{\{s<\zeta\}} \right]
=\E_x\left[\int^t_0\int_E f(X_s, y) N(X_s, dy)dH_s\right].
\end{equation}
The Revuz measure for $H$ will be denoted as $\mu_H$.

\medskip

\begin{defn}\label{D:2.2}
Suppose $F$ is a bounded function
on $E\times E$ vanishing on the diagonal $d$.
It is always extended to be zero off $E\times E$.
Define $\mu_F (dx):= \left(
\int_E F(x, y) N(x, dy)\right) \mu_H (dx)$.
We say $F$ belongs to the class $\J(X)$ (respectively, $\J_{\infty}(X)$) if
the measure
$$\mu_{|F|} (dx):= \left( \int_{E} |F(x, y)| N(x, dy)\right)
                     \mu_H (dx)
$$
belongs to $\K(X)$ (respectively, $\K_{\infty}(X)$).
\end{defn}

\medskip

See \cite[Section 2]{C2} for concrete examples of
$\mu \in \K_\infty (X)$ and $F\in \J_\infty (X)$.

\medskip

For $\alpha>0$, let $X^{(\alpha)}$ denote the $\alpha$-subprocess of $X$;
that is,  $X^{(\alpha)}$ is the subprocess of $X$ killed at exponential rate $\alpha$.
Let $G^{(\alpha)}$ be the $0$-resolvent (or Green operator) of $X^{(\alpha)}$.
Then $G^{(\alpha)}=G_\alpha$,  the $\alpha$-resolvent of $X$.
Thus for $\beta >\alpha >0$,
 $ \K_1 (X)\subset \K_1 (X^{(\alpha )}) \subset
\K_1 (X^{(\beta )})$ and
$\K_\infty (X)\subset \K_\infty (X^{(\alpha )}) \subset
\K_\infty (X^{(\beta )}) $.
In fact, it follows from the resolvent equation
$G_\alpha =G_\beta + (\beta -\alpha) G_\alpha G_\beta$
that $\K_\infty (X^{(\alpha)}) = \K_\infty (X^{(\beta)})$
for every $\beta >\alpha>0$.
Consequently,  $\J_\infty (X^{(\alpha)}) = \J_\infty (X^{(\beta)})$
for every $\beta >\alpha>0$. Clearly, $\K (X^{(\alpha)})=\K (X)$
for every $\alpha>0$.

Assume that $\mu$ is a signed smooth measure with $\mu^+\in
\K  (X)$ and $G\mu^-$ bounded, and $F\in \J  (X)$ symmetric.
Define the non-local Feynman-Kac semigroup
$$P^{\mu, F}_tf (x):=\bE_x \Big[ \exp\Big( A^\mu_t+\sum_{0<s\leq t} F(X_{s-}, X_s) \Big) f(X_t)\Big], \qquad t\geq 0.
$$
It follows from the proof of \cite[Proposition 2.3]{CS1}
and H\"older  inequality that
$\{P^{\mu, F}_t; t\geq 0\}$ is a strongly continuous semigroup on
  $L^p(E; m) $ for every $1\leq p\leq \infty$. Moveover, it is easy to verify that $P^{\mu, F}_t$ is a symmetric operator in $L^2(E; m)$ for every
  $t\geq 0$.
The $L^p$-spectral bound of $\{ P^{\mu, F}_t; t\geq 0\}$ is defined to
be
$$ \lambda_p(X,  \mu+F):= -\lim_{t\to \infty}
    \frac1t \log \| P^{\mu, F}_t\|_{p, p}.
$$

Necessary and sufficient conditions
for $\lambda_p(X,  \mu+F)$ to be independent of $1\leq p\leq \infty$
have been investigated in \cite{C3} by using gaugeability results
for Schr\"odinger semigroups  obtained in \cite{C}.
The following three results are established in \cite{C3}.

\begin{thm}\label{T:2.3}
{\rm (See \cite[Theorem 5.3]{C3})}
Assume that    $m(E)<\infty$ and that the following condition holds
\begin{equation}\label{e:lp3}
 \hbox{there is some } t_0> 0 \hbox{ so that }
 P_{t_0} \hbox{ is a bounded operator from }
 L^2(E; m) \hbox{ into } L^\infty (E; m).
\end{equation}
Let $\mu$ be  a signed smooth measure
with $\mu^+\in  \K_\infty (X^{(\alpha)})$
and $G_{\alpha } \mu^-$ bounded for some $\alpha \geq 0$,
and $F\in \J_\infty (X^{(\alpha)})$ symmetric.
Then
$\lambda_p(X, \mu+F)$ is independent of $p\in [1, \infty]$.
\end{thm}

\begin{thm}\label{T:2.4}
{\rm (See \cite[Theorem 5.4]{C3})}
 Suppose that $\mu$ is a signed smooth measure with
$\mu^+\in \K_\infty (X^{(1)})$ and $G_1\mu^-$ bounded,
and $F\in \J_\infty (X^{(1)})$ symmetric.
\begin{description}
\item{\rm (i)} $\lambda_\infty (X,  \mu+F)\geq \min \{\lambda_2 (X, \mu+F), 0\}$. Consequently, $\lambda_p (X,  \mu+F)$ is independent of $p\in [1, \infty]$ if $\lambda_2 (X, \mu+F)\leq 0$.
\item{\rm (ii)} Assume in addition that  $X$ is conservative
and that   $\mu\in \K_\infty (X^{(1)})$.
 Then $\lambda_\infty (X, \mu +F )=0$
 if   $\lambda_2 (X, \mu+F)>0$.
 Hence   $\lambda_p (X, \mu +F)$ is independent of $p\in [1, \infty]$ if and only if $\lambda_2(X, \mu +F )\leq 0$.
\end{description}
\end{thm}

\begin{thm}\label{T:2.5}
{\rm (\cite[Theorem 5.5]{C3})}
 Suppose that $1\in \K_\infty (X^{(1)})$,
 $\mu  \in \K_\infty (X^{(1)})$ and $F\in \J_\infty (X^{(1)})$
 symmetric.
Then  $\lambda_p(X, \mu+F)$ is independent of $p\in [1, \infty]$.
\end{thm}

\section{Generalized Feynman-Kac semigroup}\label{S:3}

Denote by $\sF_e$ the extended Dirichlet space of $(\sE, \sF)$.
Every $u\in \sF_e$ admits a quasi-continuous version, which we still denote
as $u$. In this paper, every $u\in \sF_e$ is always represented by its
quasi-continuous version. For such $u$, the following Fukushima's decomposition holds (cf. \cite{CF, FOT}):
$$ u(X_t)=u(X_0)+M^u_t + N^u_t, \qquad t\geq 0,
$$
where $M^u$ is a martingale additive functional of $X$ having finite energy
and $N^u$ is a continuous additive functional of $X$ having zero energy.
The continuous martingale part of $M^u$ will be denoted as $M^{u, c}$.
Let $\< M^u\>$ and $\<M^{u, c}\>$ be the predictable quadratic variation
processes of $M^u$ and $M^{u, c}$, respectively. Both of them are
positive continuous additive functionals of $X$, whose Revuz measures
will be denoted as $\mu_{\<u\>}$ and $\mu_{\<u\>}^c$, respectively.
Note that by \cite[Theorem 4.3.11]{CF} or \cite[Theorem 5.2.3]{FOT}, for bounded $u$ in $\sF_e$,
$\mu_{\< u\>}$ can be  computed from
  \begin{equation}\label{e:3.1a}
  \int_E f(x) \mu_{\< u\>}(dx) = 2\sE (uf, u)-\sE (u^2, f)
  \qquad \hbox{for bounded } f\in \sF_e.
  \end{equation}
A similar formula holds for $\mu_{\<u\>}^c$ as well; see
\cite[Exercise 4.3.12]{CF}.

  Let $u$ be a bounded function in $\sF_e$ with
  $\mu_{\< u\>}\in \K_\infty (X^{(1)})$, $\mu \in \K_\infty (X^{(1)})$
  and $F$ be a bounded symmetric function  in $\J_\infty (X^{(1)})$.
  Define the Feynman-Kac semigroup $\{T^{u, \mu, F}_t, t\geq 0\}$ by
$$
  T^{u, \mu, F}_t f(x)=\E_x \left[ \exp \left( N^u+A^\mu_t+\sum_{0<s\leq t}
   F(X_{s-}, X_s)\right) f(X_t) \right] .
$$
We will show that for every $p\in [1, \infty]$,
$\{T^{u, \mu, F}_t, t\geq 0\}$
is a strongly continuous symmetric semigroup on $L^2(E; m)$.
This will be achieved by reducing the generalized non-local
Feynman-Kac semigroup $\{T^{u, \mu, F}_t, t\geq 0\}$
via a suitable Girsanov transform to a non-local
Feynman-Kac semigroup of the new process.

Note that since $u$ is bounded, $v (x):= e^{-u}-1$ is a bounded function in
$\sF_e$. Clearly for every $x, y\in E$,
$$ |v (x) | \leq e^{\| u\|_\infty} \, | u(x) |
\quad \hbox{ and so } \quad |v(x)-v(y)|\leq e^{\| u\|_\infty} \, | u(x) -u(y)|.
$$
We thus deduce from \cite[(4.3.12) and Theorem 4.3.7]{CF} that
\begin{equation}\label{e:3.1}
\mu_{\< v \>} (dx) \leq e^{2 \| u\|_\infty} \, \mu_{\< u\>} (dx).
\end{equation}
Let $Z={\rm Exp} (M)$ be the Dol\'eans-Dade exponential martingale of $M_t:=\int_0^t e^{u(X_{s-})} dM^v_s$; that is, $Z$ is the unique
solution of
$$ Z_t=1+\int_0^t Z_{s-}   dM_s, \qquad t\geq 0.
$$
It follows from Dol\'eans-Dade formula (cf. \cite[Theorem 9.39]{HWY})
that
\begin{eqnarray}
Z_t&=& \exp\left( M_t-{1\over 2} \< M^c\>_t \right)
\prod_{0<s\leq t} \left( 1+ M_s- M_{s-}\right)
e^{-( M_s- M_{s-})} \nonumber \\
&=& \exp\left(M_t-{1\over 2}\<M^{u,c}\>_t \right)
\prod_{0<s\leq t}
\exp \left(u(X_{s-})-u(X_{s})+1-  e^{u (X_{s-})-u (X_{s})} \right).
\label{e:3.2}
\end{eqnarray}
Note that
$$ M_t-M_{t-}= e^{u(X_{t-})-u(X_t)}-1\geq e^{-2\| u \|_\infty} -1$$
and
that by \eqref{e:3.1},
$$ \sup_{x\in E} \E_x [ M ]_\infty
=  \sup_{x\in E} \E_x \left[ \int_0^\infty e^{2u(X_{s-})} d [M^v]_s\right]
\leq e^{2\| u\|_\infty} \sup_{x\in E} \E_x \left[\< M^v \>_\infty \right]
\leq e^{4\| u\|_\infty} \| G \mu_{\< u\>} \|_\infty <\infty,
$$
where $[M]$ is the quadratic variation process of the martingale $M$.
Therefore we conclude by the uniform integrability criteria for exponential
 martingales established in \cite[Theorem 3.2]{C3} that
$Z={\rm Exp}(M)$ is a uniformly integrable martingale under $\P_x$
for every $x\in E$.

Let $\{\wt \P_x, x\in E\}$ be the family of probability measures defined by
$$ \frac{d\wt \P_x}{d\P_x} = Z_\infty \qquad  \hbox{on } \sF_\infty,
$$
and for emphasis, let $\wt X=(\wt X_t, \wt \P_x)$ denote the Girsanov transformed
process $(X_t, \wt \P_x)$. The following  result is  proved in \cite[Theorem 3.4]{CZ}.

\begin{thm} The process $\wt X$ is a symmetric strong Markov process with symmetrizing
measure $e^{-2u(x)} m(dx)$, whose associated Dirichlet form on $L^2(E; e^{-2u(x)}m(dx))$
 is $(\wt  \sE, \sF)$, where for $f\in \sF$,
\begin{eqnarray*}
\wt \sE (f,f) &=&
\frac12 \int_E  e^{-2u(x)} \mu^c_{\<f\>}(dx)
+\frac12 \int_{E\times E\setminus d}
 ( f(x)-f(y))^2  e^{ -u (x) -u (y)} N(x, dy) \mu_H (dx) \\
 &&  +\int_E  f(x)^2 e^{-u(x)} \kappa (dx) .
\end{eqnarray*}
\end{thm}

It follows that $\wt X$ has a L\'evy system $(\wt N(x, dy), \wt H)$ with
$$
\wt N(x, dy) = e^{-u(y)} N(x, dy) \quad \hbox{ and }
\quad \mu_{\wt H} (dx)= e^{-u(x)} \mu_H (dx).
$$
 In view of Lemma \ref{L:3.2} below, the latter is equivalent to
$\wt H_t= \int_0^t e^{u(X_s)} dH_s$.
The next lemma is established in \cite[Theorem 3.3 and Lemma 4.2]{CZ}.

\begin{lemma}\label{L:3.2}
 If $A$ is a positive continuous additive functional of $X$ with
Revuz measure $\mu$, then $A$ is a  positive continuous additive functional of $\wt X$
 with Revuz measure $e^{-2u(x)}\mu (dx)$. Moreover, if $\mu \in \K  (X)$, then
 $e^{-2u(x)}\mu (dx) \in \K  (\wt X )$.
\end{lemma}

Since $u$ is bounded, the second half of the above lemma says that
  $\K  (X) \subset \K  (\wt X )$.

\begin{lemma}\label{L:3.3}
$$
\K_\infty (X^{(1)}) \subset \K_\infty (\wt X^{(1)})
\quad \hbox{and so } \quad
\J_\infty (X^{(1)}) \subset \J_\infty (\wt X^{(1)}).
$$
\end{lemma}

\pf Let $\mu$ be a non-negative measure in $\K_\infty (X^{(1)})$.
It suffices to show that $\nu (dx):=e^{-2u(x)} \mu (dx)
\in \K_\infty (\wt X^{(1)})$.
Let $A$ be the positive  continuous additive functional of $X$ having Revuz measure $\mu$. In view of Lemma \ref{L:3.2}, it
can also be viewed as the positive continuous additive
functional of $\wt X$ with Revuz measure $\nu$.
Observe that $\mu \in \K (X^{(1)})= \K (X)$.
By Lemma \ref{L:3.2}, $\nu \in \K (\wt X)$ and
so there is $\alpha>0$ such that
$\| \wt G_\alpha \nu \|_\infty \leq 1$.
For any given $\eps>0$, choose $t_0>0$ so that $e^{-\alpha t_0} <\eps /2$.
Then for any $x\in E$ and any $B\in \sB(E)$,
\begin{eqnarray}
\wt G_\alpha (\1_B \nu) (x)
&=& \E_x \left[ \int_0^\infty e^{-\alpha s} \1_B(\wt X_s) dA_s \right]
\leq   \wt \E_x \left[ \int_0^{t_0} \1_B (\wt X_s) dA_s \right]
+ \wt \E_x \left[\int_{t_0}^\infty e^{-\alpha s} dA_s\right] \nonumber \\
&\leq &  \E_x \left[ Z_t \int_0^{t_0} \1_B(X_s) dA_s \right] +
e^{-\alpha t_0} \wt \E_x \left[
 \wt G_\alpha
\mu (X_{t_0}) \right] \nonumber \\
&\leq & \left( \E_x \left[ Z_{t_0}^2\right]\right)^{1/2}  \left(
\E_x \left[ \left( \int_0^{t_0} \1_B(X_s) dA_s \right)^2\right] \right)^{1/2}
+ \eps /2. \label{e:3.2b}
\end{eqnarray}
By \cite[Lemma 4.1(ii)]{CZ}, $\sup_{x\in E} \E_x [ Z_{t_0}^2]=c(t_0) < \infty$.
On the other hand, denoting $f(x):=\E_x  \left[\int_0^{t_0} \1_B(X_s) dA_s\right]$,
we have by the Markov property of $X$,
\begin{eqnarray*}
&& \E_x \left[ \left( \int_0^{t_0} \1_B(X_s) dA_s \right)^2\right]
=  2 \E_x \left[ \int_0^{t_0} \1_B(X_s) \left( \int_s^{t_0} \1_B(X_r) dA_r\right) dA_s\right] \\
&\leq & 2\E_x \left[ \int_0^{t_0} \1_B(X_s) f(X_s) dA_s\right]
\leq 2 \| f\|_\infty^2 \leq 2 e^{2\alpha t_0} \| G_\alpha (\1_B \mu) \|_\infty^2.
\end{eqnarray*}
This together with \eqref{e:3.2b} yields
$$
\| \wt G_\alpha (\1_B \nu)\|_\infty \leq \sqrt{2c(t_0)}\, e^{\alpha t_0}
 \, \| G_\alpha (\1_B \mu) \|_\infty +\eps/2.
$$
It follows from the definition of $\K_\infty$ that $\nu \in \K_\infty (\wt X^{(\alpha)})
=   \K_\infty (\wt X^{(1)})$. \qed

\begin{remark} \rm
(i) By \cite[Theorem 3.5]{CZ}, $X$ can be recovered from $\wt X$ through an analogous Girsanov
transform. Thus we in fact have  $\K  (X) = \K  (\wt X )$,
 $\K_\infty (X^{(1)}) = \K_\infty (\wt X^{(1)})$ and
$\J_\infty (X^{(1)}) = \J_\infty (\wt X^{(1)})$.

See \cite[Lemma 3.3]{DKK} for a related result
 under the assumption that $X$ is a doubly Feller process
with no killings inside.
\qed
\end{remark}

In view of \eqref{e:3.2}, we can express $e^{N^u_t}$ as follows (see \cite[(4.6)]{CZ}),
\begin{equation}\label{e:3.4}
 \exp (N^u_t)= \exp \left( u(X_t)-u(X_0)-M^u_t \right)=
 e^{-u(x)}  Z_t \, e^{-A_t} \, e^{u(X_t)}  ,
\end{equation}
where
 $$
A_t: = \int_0^t\left(\int_{E_\partial}
\left(  u(X_s)- u(y) +1-e^{u(X_s) -  u(y)  } \right)
N(X_s,dy)\right)dH_s  -{1\over 2}\<M^{u,c}\>_t .
$$
Let $\nu$ be the signed Revuz measure of $A$, that is,
\begin{equation}\label{e:3.6a}
\nu (dx):= \left(\int_{E_\partial}
\left(  u(x)- u(y) +1-e^{u(x) -  u(y)  } \right)
N(x,dy)\right)d\mu_H  (dx) -{1\over 2} \mu^c_{\< u \>} (dx).
\end{equation}
Since $u$ is bounded,
\begin{equation}\label{e:3.5}
| \nu (dx)| \leq \frac{e^{\| u\|_\infty}}2 \left(\int_{E_\partial}
(u(x)-u(y))^2 N(x, dy)\right) d\mu_H (dx) +
\frac12 \mu^c_{\< u \>} (dx)\leq \frac{e^{\| u\|_\infty}}2 \mu_{\< u\>}(dx)
\end{equation}
and so $\nu \in \K_\infty (X^{(1)})\subset \K_\infty (\wt X^{(1)})$.

For convenience, if $A^\mu$ is a continuous additive functional of $X$ with
(signed) Revuz measure $\mu$, in view of Lemma \ref{L:3.2},
we will denote $A^\mu$ by $\wt A^{\, e^{-2u}\mu}$
when viewed as a continuous additive functional of $\wt X$.

By \eqref{e:3.4},
\begin{eqnarray}
T^{u, \mu, F}_t f(x) &=&  e^{-u(x)} \E_x\left[   Z_\infty \,
\exp \left(-A^\nu_t + A^\mu_t + \sum_{0<s\leq t}
F (X_{s-}, X_s) \right)\, (e^{u} f)(X_t) \right] \nonumber \\
&=&  e^{-u(x)} \wt \E_x\left[
\exp \left( \wt A^{\, e^{-2u}(\mu -\nu)}_t + \sum_{0<s\leq t}
F (\wt X_{s-}, \wt X_s) \right)\, (e^{u} f)(\wt X_t) \right]
\nonumber \\
&=& e^{-u(x)}\, \wt T_t^{\, e^{-2u}(\mu -\nu), F} (e^u f)(x),
\label{e:3.6}
\end{eqnarray}
 where $\{\wt T_t^{\, e^{-2u}(\mu -\nu), F}, t\geq 0\}$
 is the non-local Feynman-Kac semigroup of $\wt X$ defined by
 $$  \wt T_t^{\, e^{-2u}(\mu -\nu), F} g(x)=\wt \E_x \left[
  \exp \left( \wt A^{\, e^{-2u}(\mu -\nu)}_t + \sum_{0<s\leq t}
    F (\wt X_{s-}, \wt X_s) \right)\, g (\wt X_t) \right] .
 $$

\begin{thm}\label{T:3.5}
  Let $u$ be a bounded function in $\sF_e$ with
  $\mu_{\< u\>}\in \K_\infty (X^{(1)})$, $\mu \in \K_\infty (X^{(1)})$
  and $F$ be a bounded symmetric function  in $\J_\infty (X^{(1)})$.
Then for every $p\in [1, \infty]$,  $\{T^{u, \mu, F}_t, t\geq 0\}$
is a strongly continuous symmetric semigroup on $L^p(E; m)$.
\end{thm}

\pf Since by \eqref{e:3.5} and Lemma \ref{L:3.3},
$\nu\in \K_\infty (\wt X^{(1)})$, $\mu\in \K_\infty (  X^{(1)})\subset
\K_\infty (\wt X^{(1)})$ and $F\in \J_\infty (\wt X^{(1)})$,
it follows from the proof of \cite[Proposition 2.3]{CS1}
 applied to the process $\wt X$ that
for every $p\in [1, \infty]$,
$\big\{ \wt T_t^{\, e^{-2u}(\mu -\nu), F}, t\geq 0 \big\}$
 is a strongly continuous symmetric semigroup on $L^p(E; e^{2u} dm)$.
Thus we have by \eqref{e:3.6} that   $\{ T_t^{u, \mu, F}, t\geq 0\}$
 is a strongly continuous symmetric semigroup on $L^p(E;   dm)$ for every
 $p\in [1, \infty]$.
 \qed

Denote the operator norm of
$\wt T_t^{\, e^{-2u}(\mu -\nu), F}: L^p(E; e^{2u} dm)
\to L^p(E; e^{2u} dm)$ by   $\big\| \wt T_t^{\, e^{-2u}(\mu -\nu), F}
\big\|_{p, p}$,
and  the operator norm of $T_t^{u, \mu, F}: L^p(E; m)\to L^p(E; m)$ by
  $\|   T_t^{u, \mu, F}\|_{p, p}$.
For $1\leq p \leq \infty$,  the $L^p$-spectral
bound of semigroup $T^{u, \mu, F}_t$ is defined as
$$ \lambda_p(X, u+\mu+F):= - \lim_{t\to \infty} \frac1t
\log \|T^{u, \mu, F} _t \|_{p,p}.
$$
Clearly, in view of \eqref{e:3.6},
 $$ \|  T_t^{u, \mu, F}\|_{2, 2} = \| \wt T_t^{e^{-2u}(\mu -\nu), F}\|_{2, 2},
 $$
 while
 $$
  e^{-2\| u\|_\infty} \|  T_t^{u, \mu, F}\|_{\infty, \infty}
   \leq  \| \wt T_t^{e^{-2u}(\mu -\nu), F}\|_{\infty, \infty}
   \leq  e^{2\| u\|_\infty} \|  T_t^{u, \mu, F}\|_{\infty, \infty}.
 $$
 It follows that
 \begin{equation}\label{e:3.8}
  \lambda_2 (X, u+\mu+F)= \lambda_2 (\wt X, e^{-2u}(\mu -\nu)+ F)
 \quad \hbox{and} \quad
  \lambda_\infty (X, u+\mu+F)= \lambda_\infty (\wt X, e^{-2u}(\mu -\nu)+ F).
 \end{equation}
Note that by    \cite[(5.7)]{C3}
\begin{eqnarray}
 &&\lambda_2 (X, u+\mu+F)= \lambda_2 (\wt X, e^{-2u}(\mu -\nu), F) \nonumber \\
&=& \inf \left\{ \wt \sE (g, g) -
 \int_{E\times E} g(x)g(y) \left( e^{F(x, y)}-1\right)e^{-u(x)-u(y)} N(x, dy) \mu_H(dx) \right. \nonumber \\
&&  \hskip 0.5truein \left.  -\int_E g(x)^2 e^{-2u(x)} (\nu -\mu) (dx); \
 g\in \sF \ \hbox{ with } \int_E g(x)^2 e^{-2u(x)} m(dx)=1 \right\}  \nonumber \\
&=& \inf \left\{ \wt \sE (g, g) -
 \int_{E\times E} g(x)g(y)e^{-u(x)-u(y)} \left( e^{F(x, y)}-1\right) N(x, dy) \mu_H(dx) \right. \nonumber \\
&&  \hskip 0.5truein \left.  -\int_E g(x)^2 e^{-2u(x)} (\nu -\mu) (dx); \
 g\in \sF_b \ \hbox{ with } \int_E g(x)^2 e^{-2u(x)} m(dx)=1 \right\}  \nonumber \\
&=& \inf \left\{ \wt \sE (fe^u, fe^u) -
 \int_{E\times E} f(x)f(y) \left( e^{F(x, y)}-1\right)  N(x, dy) \mu_H(dx) \right. \nonumber \\
&&  \hskip 0.5truein \left.  -\int_E f(x)^2   (\nu -\mu) (dx); \
 f\in \sF_b \ \hbox{ with } \int_E f(x)^2   m(dx)=1 \right\}  \nonumber \\
 &=& \inf \left\{   \sE (f, f) + \sE(u, f^2) +\int_E f(x)^2   \mu (dx)
 -\int_{E\times E} f(x)f(y) \left( e^{F(x, y)}-1\right)  N(x, dy) \mu_H(dx) ; \right. \nonumber \\
&&  \hskip 0.8truein \left.  \
 f\in \sF_b \ \hbox{ with } \int_E f(x)^2   m(dx)=1 \right\}  .
 \label{e:3.9}
\end{eqnarray}
In the last equality, we used the fact that
$$ \sE (f e^u, fe^u) -\int_E f(x)^2 \nu (dx) =
\sE(f, f)+\sE (u, f^2) \qquad \hbox{for } f\in \sF_b,
$$
whose proof can be found in the paragraph following
(4.8) in the proof of \cite[Theorem 1.2]{CZ} for bounded $u$.

Clearly
$$ \| T^{u, \mu, F}_t \|_{\infty, \infty} =
\| T^{u, \mu, F}_t 1\|_\infty =\sup_{x\in E}
\E_x \left[ \exp \left( N^u+A^\mu_t+\sum_{0<s\leq t}
   F(X_{s-}, X_s)\right); t<\zeta  \right].
$$
By duality, we have  $\| T^{u, \mu, F}_t\|_{1,1}= \| T^{u, \mu, F}_t \|_{\infty, \infty}$.
Consequently, it follows from the Cauchy-Schwarz inequality that
$$
 \|T^{u, \mu, F}_t f\|^2_2 \leq  \| T^{u, \mu, F}_t 1\|_\infty
  \, \| T^{u, \mu, F}_t (f^2)\|_1
\leq \| T^{u, \mu, F}_t\|_{\infty, \infty}^2
\| f\|^2_2 \qquad \hbox{for } f\in L^2(E; m).
$$
Thus we have $\| T^{u, \mu, F}_t \|_{2, 2}  \leq
\|T^{u, \mu, F}_t\|_{\infty, \infty}$.
We now deduce by interpolation that
$$
\| T^{u, \mu, F}_t \|_{2,2}\leq \| T^{u, \mu, F}_t\|_{p,p} \leq
\| T^{u, \mu, F}_t \|_{\infty, \infty}
\quad \mbox{for } 1< p< \infty.
$$
Hence
\begin{equation}\label{e:4.2}
\lambda_\infty (X, u+\mu+F) \leq \lambda_p (X, u+\mu+F)
\leq \lambda_2 (X, u+\mu +F) \quad \mbox{for } 1<p<\infty.
\end{equation}

\section{$L^p$-independence of spectral bounds}\label{S:4}

We can now present results on the $L^p$-independence of
the spectral bounds of generalized non-local Feynman-Kac
semigroups.

\begin{thm}\label{T:eigen}
Assume that   $m(E)<\infty$,   and that \eqref{e:lp3} holds.
Let $u$ be a bounded function in $\sF_e$ with
  $\mu_{\< u\>}\in \K_\infty (X^{(1)})$,
   $\mu\in \K_\infty (X^{(1)})$ and $F$ a symmetric
function in $\J_\infty (X^{(1)})$.
Then
$\lambda_p(X, u+\mu+F)$ is independent of $p\in [1, \infty]$.
\end{thm}

\pf
Since $P_{t_0}$ is a bounded linear operator from $L^2(E; m)$
  to $L^\infty (E; m)$, by duality, $P_{t_0}$ is a bounded linear operator
  from $L^1(E; m)$ to $L^2(E; m)$. Hence $P_{2t_0}: L^1(E; m)\to L^\infty (E; m)$ is bounded. Let $M_t:=\int_0^t e^{u(X_{s-})} dM^{e^{-u}-1}_s$
  and $Z_t =\Exp (M)_t$ be its Dol\'eans-Dade exponential martingale,
   which admits expression \eqref{e:3.2}.
  Since $\mu_{\<u \>}\in \K_\infty (X^{(1)})\subset\K (X^{(1)})= \K (X)$,
   we have by \cite[Lemma 4.1(ii)]{CZ} that
$\sup_{x\in E} \E \left[ Z_{2t_0}^2\right] <\infty$.

Denote by $Y$ the Girsanov transformed process of $X$ via $Z$.
 Then for every $f\in L^2(E; m)$,
$$ |P^Y_{2t_0} f(x)|:= |\E_x [ f(Y_{2t_0})]|=
\left| \E_x \left[ Z_{2t_0} f(X_t)\right] \right|
\leq \left( \E_x [Z^2_{2t_0}]
 \, \E_x \left[ f(X_{2t_0})^2\right]\right)^{1/2}
 \leq c\, \|f\|_{L^2(E; m)}.
 $$
 This proves that condition \eqref{e:lp3} holds for $Y$
  with $2t_0$ in place of $t_0$. Let $\nu$ be the measure defined
  by \eqref{e:3.6a}.  Note that in view of \eqref{e:3.5}, $\nu\in  \K_\infty (X^{(1)})\subset \K_\infty (\wt X^{(1)})$.
  Thus we can apply Theorem \ref{T:2.3} to conclude that
  $$
  \lambda_2(Y, e^{-2u}(\mu -\nu)+ F)
   =\lambda_\infty (Y, e^{-2u}(\mu -\nu)+ F).
  $$
  We deduce from this and \eqref{e:3.8} that
  $\lambda_2(X, u+\mu+ F)
  =\lambda_\infty (u+\mu+ F)$.
  Consequently, we have from \eqref{e:4.2} that
  $\lambda_p(X, u+\mu+ F)$ is independent of $p\in [1, \infty]$.
     \qed

\begin{thm}\label{T:5.4}
 Suppose that   $u$ is a bounded function in $\sF_e$ with
  $\mu_{\< u\>}\in \K_\infty (X^{(1)})$, $\mu \in \K_\infty (X^{(1)})$
  and $F\in \J_\infty (X^{(1)})$ is bounded and symmetric.
\begin{description}
\item{\rm (i)} $\lambda_\infty (X,  u+\mu+F)\geq \min \{\lambda_2 (X, u+ \mu+F), 0\}$ and so $\lambda_p (X, u+\mu+F)$ is independent of $p\in [1, \infty]$ if $\lambda_2 (X, u+\mu+F)\leq 0$.
\item{\rm (ii)} Assume in addition that  $X$ is conservative.
 If $\lambda_2 (X, u+\mu+F)>0$, then $\lambda_\infty (X, u+\mu +F )=0$.
 Hence   $\lambda_p (X, u+\mu +F)$ is independent of $p\in [1, \infty]$ if and only if $\lambda_2(X, u+ \mu +F )\leq 0$.
\end{description}
\end{thm}

\bigskip

\pf   Let $Z_t=\Exp (M)_t$ be the exponential martingale in
the proof of Theorem \ref{T:eigen}.
As we saw in Section \ref{S:3},
 $\{Z_t, t\geq 0\}$ is a uniformly integrable martingale under each $\bP_x$. It follows that the Girsanov transformed process $\wt X$ of $X$ by $Z$
 is transient and has a Green function
$G^Y$. Furthermore, $\wt X$ is conservative if so is $X$.
  It is clear that $\wt X$ is $e^{-2u}m$-irreducible since $\Exp(M)_t>0$ a.s..
Let $\nu$ be the measure defined by \eqref{e:3.6a}.
Note that in view of \eqref{e:3.5}, $\nu\in  \K_\infty (X^{(1)})\subset \K_\infty (\wt X^{(1)})$.
 Since $e^{-2u}(\mu-\nu)\in \K_\infty (X^{(1)})\subset \K_\infty (\wt X^{(1)})$,
 the conclusion of the theorem now follows
from  Theorem \ref{T:2.4} applied to $(\wt X , e^{-2u}(\mu-\nu)+ F)$
and \eqref{e:3.8}. \qed

\begin{thm}\label{T:5.5}
 Suppose that $1\in \K_\infty (X^{(1)})$,  $u$ is a bounded function
 in $\sF_e$ with $\mu_{\< u\>}\in \K_\infty (X^{(1)})$,
 $\mu  \in \K_\infty (X^{(1)})$ and $F\in \J_\infty (X^{(1)})$
 symmetric.
Then  $\lambda_p(X, u+\mu+F)$ is independent of $p\in [1, \infty]$.
\end{thm}

\pf  As above, let $Z_t=\Exp (M)_t$ be the exponential martingale in
the proof of Theorem \ref{T:eigen} and $\wt X$ the Girsanov
 transformed process of $X$ by $Z$. Let $\nu$ be the measure defined
 by \eqref{e:3.6a}, which in view of  \eqref{e:3.5} is in
  $\K_\infty (X^{(1)})\subset \K_\infty (\wt X^{(1)})$.
 The conclusion of this theorem follows from  Theorem \ref{T:2.5} applied to $(\wt X, e^{-2u}(\mu-\nu)+ F)$ and \eqref{e:3.8}. \qed

\section{Examples}\label{S:5}

In this section, we give several concrete examples for functions
to be in Kato class $\K_\infty (X^{(1)})$, $\J_\infty (X^{(1)})$
and for bounded $u\in \sF_e$ with $\mu_{\<u\>}\in \K_\infty (X^{(1)})$
so that the main results of this paper apply.

Two real-valued functions $f$ and $g$ are said to be
comparable if there is a constant $c>1$ so that
$g/c\leq f \leq c g$, and we denote it by $f\asymp g$.

\bigskip

\noindent{\bf Example 5.1} {\bf (Stable-like process on $d$-sets)}
Let $n\geq 1$ and $0<d\leq n$.
A Borel subset  $E$ in $\bR^n$   is said to be a global $d$-set if
 there exist a measure $m$ on $E$ and
  constants   $C_2> C_1>0$ so that
\begin{equation}\label{eqn:dset}
C_1\, r^d\le m (B(x,r) )\le C_2 \, r^d \qquad \hbox{for all }
 x\in E \hbox{ and } r>0.
\end{equation}
Here $B(x,r):=\{y\in E: |x-y|<r\}$ and
$|\cdot|$ is the Euclidean metric in $\bR^n$.

For a closed global $d$-set $E\subset \bR^n$ and $0<\alpha <2 $,
 define
\begin{eqnarray}
\sF  &=& \left\{ u \in L^2(E, m): \, \int_{E\times E}
\frac{(u(x)-u(y))^2}
{|x-y|^{d+\alpha}} \, m (dx) m(dy) < \infty \right\}
\label{eqn:form1} \\
\sE(u,v)&=& \frac 12 \int_{E\times E} (u(x)-u(y))(v(x)-v(y)) \frac{
c(x, y) } {|x-y|^{d+\alpha}} \, m(dx) m (dy)  \label{eqn:form2}
\end{eqnarray}
for $u, \, v\in \sF$,
where $c(x, y)$ is a symmetric function on $E\times E$ that is bounded
between two strictly positive constants $C_4>C_3>0$, that is,
\begin{equation}\label{eqn:1.4}
C_3 \leq c(x, y) \leq C_4 \qquad \hbox{for } m \hbox{-a.e. }
   x, y \in E.
\end{equation}
It is shown in \cite{CK1} that   $(\sE, \sF)$ is a regular Dirichlet form on
$L^2(E; m)$ and  there is an associated $m$-symmetric Hunt
process $X$ on $E$ starting from every point in $E$.
Moreover, $X$ admits a jointly H\"older continuous transition
density function $p(t, x, y)$ with respect to the measure $m$,
which satisfies the following two-sided estimates
\begin{equation}\label{e:5.5}
p(t, x, y) \asymp t^{-d/\alpha} \wedge \frac{t}{|x-y|^{d+\alpha}}
\qquad \hbox{on } (0, \infty) \times E \times E,
\end{equation}
where the comparison constants in \eqref{e:5.5} depend only on
$C_k$, $k=1, 2, 3, 4$.
We call such kind of process
a $\alpha$-stable-like process on $E$.
 Note that when $E=\bR^n$ and
$c(x, y)$ is a constant function, then $X$ is nothing but a
rotationally symmetric $\alpha$-stable process on $\bR^n$.
The process $X$  has a L\'evy system $(N(x, dy), H)$, with
$N(x, dy)=N(x, y) dy=c(x, y)
 |x-y|^{-(n+\alpha)} dy$  and $H_t=t$ so $\mu_H(dx)=m (dx)$.
When $\alpha <d$, the process $X$ is transient and its Green function
is given by
$$
  G(x, y) =\int_0^\infty p(t, x, y) dt \asymp
 |x-y|^{\alpha -d} , \qquad x, y \in E.
$$
By the same argument as that for Theorem 2.1 of Chen \cite{C},
we can show that when $0<\alpha < d\wedge 2$,
\begin{description}
\item{(a)}
a  signed measure $\mu$ is in $\K (X)$
if and only if
\begin{equation}\label{eqn:2.1}
 \lim_{r\to 0} \sup_{x\in E} \int_{B(x,r)}
|x-y|^{\alpha -d} |\mu|(dy)=0 ;
\end{equation}

\item{(b)} a finite signed measure $\mu$ is in
$\K_\infty (X)$ if and only if it is in $\K(X)$;

\item{(c)}
a signed measure $\mu$ is in $\K_\infty (X)$
if and only if   both (\ref{eqn:2.1}) and the following
condition
\begin{equation}\label{eqn:2.2}
 \lim_{R\to \infty} \sup_{x\in E} \int_{B(0,R)^c}
|x-y|^{\alpha -d} |\mu|(dy)=0
\end{equation}
are satisfied.
\end{description}
It is easy to see that condition (\ref{eqn:2.1})
is satisfied for $\mu(dx)=f(x)m (dx) $ if $f\in L^p(E; m)$ for
some $p>d/\alpha$.
We next show the following.

\begin{lemma}\label{L:5.1} Let $0<\alpha <2$ and $X$ be a symmetric
$\alpha$-stable-like process on the $d$-set $E$.
\begin{description}
\item{\rm (i)} When  $\alpha \leq d$,
$L^p(E; m)\subset  \K_\infty (X^{(1)})$ for every $p>d/\alpha$.
 When $0<d< \alpha$,
$L^p(E; m)\subset \K_\infty (X^{(1)})$ for every $p\geq 1$.

\item{\rm (ii)} For bounded $u\in \sF_e$, $\mu_{\<u\>} \in \K_\infty (X^{(1)})$
 if $f_u (x):= \int_{E} ( u(x)-u(y))^2 \frac{c(x, y)}{|x-y|^{d+\alpha}}
 m(dy)$ is in $L^p (E; m)$ for some $p>d/\alpha$.
In particular, if $u\in C^1_c(E)$,
then $\mu_{\<u\>} \in \K_\infty (X^{(1)})$.

\item{\rm (iii)} If $F$ is a bounded function on $E\times E$ with
\begin{equation}\label{e:5.8a}
|F(x, y)| \leq c\, |x-y|^\gamma \  \hbox{  for } \ x, y \in E
\quad \hbox{and} \quad
F(x, y) =0 \ \hbox{ for } (x, y)\in E\times K^c,
\end{equation}
 where $K$ is a compact subset of $E$,
 $c$ and $\gamma$ are two
positive constants such that $\gamma > \alpha$,
then $F\in \J_\infty (X^{(1)})$.
\end{description}
\end{lemma}

\pf (i) In view of \eqref{e:5.5},
\begin{equation}
G_1(x, y)= \int_0^\infty e^{-t}  p(t, x, y) dt
 \asymp  \int_0^{|x-y|^\alpha} e^{-t} \frac{t}{|x-y|^{d+\alpha}} dt
+ \int_{|x-y|^\alpha}^\infty e^{-t} t^{-d/\alpha} dt. \label{e:5.9}
\end{equation}
Observe that
\begin{equation}\label{e:5.10}
\int_0^{|x-y|^\alpha} e^{-t} \frac{t}{|x-y|^{d+\alpha}} dt
=   \frac{1-(1+|x-y|^\alpha)e^{-|x-y|^\alpha}}{|x-y|^{d+\alpha}}
\leq c_1 \frac{e^{-|x-y|^\alpha}}{|x-y|^{d-\alpha}},
\end{equation}
while
\begin{equation}\label{e:5.11}
\int_{|x-y|^\alpha}^\infty e^{-t} t^{-d/\alpha} dt
\leq \frac{\alpha}{d-\alpha} \frac{e^{-|x-y|^\alpha}}{|x-y|^{d-\alpha}}
 \qquad  \hbox{when } d>\alpha .
\end{equation}
When $d\leq \alpha$ and $|x-y|<1$,
\begin{eqnarray}\label{e:5.12}
  \int_{|x-y|^\alpha}^\infty e^{-t} t^{-d/\alpha} dt
  &\leq& \int_{|x-y|^\alpha}^1 t^{-d/\alpha} dt  +
  \int_1^\infty e^{-t}dt \nonumber \\
 & =&\begin{cases}
 \alpha  \log (1/|x-y|)+1 \qquad &\hbox{if } d=\alpha, \cr
 \frac{\alpha}{\alpha -d} \left(1-|x-y|^{\alpha -d}\right)+1
 &\hbox{if } d<\alpha,
 \end{cases}
\end{eqnarray}
while for $|x-y|\geq 1$,
\begin{equation}\label{e:5.13}
  \int_{|x-y|^\alpha}^\infty e^{-t} t^{-d/\alpha} dt
  \leq |x-y|^{-d} \int_{|x-y|^\alpha}^\infty e^{-t}dt = |x-y|^{-d}
  e^{-|x-y|^\alpha}.
\end{equation}
Thus we have by \eqref{e:5.9}-\eqref{e:5.13}
\begin{eqnarray}\label{e:5.14}
G_1(x, y)
&\leq& \begin{cases}
 c_2 \frac{e^{-|x-y|^\alpha}}{|x-y|^{d-\alpha}}
 \quad &\hbox{when } d>\alpha, \cr
c_2 \left( \log (1/|x-y|) \1_{\{|x-y|<1/2\}} + \frac{e^{-|x-y|^\alpha}}
{|x-y|^d} \1_{\{|x-y|\geq 1/2\}}\right) &\hbox{when } d=\alpha, \cr
c_2 \frac{e^{|x-y|^{-\alpha} }}{1+|x-y|^d } &\hbox{when } d<\alpha.
 \end{cases}
\end{eqnarray}
Suppose that $0<\alpha \leq d$.
For $f\in L^p(E; m)$ with $p>d/\alpha$, let $q>1$ be the conjugate of $p$,
that is, $q=p/(p-1)$. Note that $q<\frac{d}{d-\alpha}$.
We have by H\"older's inequality that
\begin{eqnarray*}
 \sup_{x\in E} \int_{B(0,R)^c} G_1(x, y) |f(y)| m (dy)
&\leq & \left( \sup_{x\in E} \int_{E} G_1(x, y)^q m(dy) \right)^{1/q}
\left( \int_{B(0, R)^c} |f(y)|^p m(dy) \right)^{1/p} \\
&\leq&  c \, \left( \int_{B(0, R)^c} |f(y)|^p m(dy) \right)^{1/p} .
\end{eqnarray*}
Hence for every $\eps >0$, there is some $R>0$ so that
$\sup_{x\in E} \int_{B(0,R)^c} G_1(x, y) |f(y)| m (dy)<\eps/2$.
On the other hand, there is $\delta >0$ so that for every
Borel set $B$ with $m(B)<\delta$,
$\left( \int_B |f(y)|^p m(dy)\right)^{1/p} < \frac{\eps}{2c} $
and so by the same argument as above,
$$  \sup_{x\in E} \int_{B} G_1(x, y) |f(y)| m (dy) < \eps /2.
$$
This shows that $f \in \K_\infty (X^{(1)})$.
Now assume that $0<d< \alpha <2$.
Using H\"older's inequality, we can deduce from above
that $L^p(E; m) \subset \K_\infty (X^{(1)})$ for every $p>1$.
We next show that $L^1((E; m) \subset \K_\infty (X^{(1)})$.
Note that since $0<d<\alpha<2$,  we have by
\eqref{e:5.14} that $G_1(x, y)$ is bounded. This in particular implies
that for $f\in L^1(E; m)$,
$$
\lim_{R\to \infty} \sup_{x\in E} \int_{B(0,R)^c} G_1(x, y) |f(y)| m (dy)
\leq \lim_{R\to \infty}  \int_{B(0,R)^c} |f(y)| m (dy) =0,
$$
and
$$
\lim_{\delta\to 0} \sup_{x\in E} \sup_{B: m(B)<\delta}
\int_{B} G_1(x, y) |f(y)| m (dy) \leq \lim_{\delta\to 0} \sup_{B: m(B)<\delta}  \int_{B(0,R)^c} |f(y)| m (dy) =0.
$$
Therefore we have $L^1(E; m)\subset \K_\infty (X^{(1)})$.

(ii) For bounded $u\in \sF_e$, we deduce from \eqref{e:3.1a} that
$$
\mu_{\< u\>} (dx) = \left( \int_{E} ( u(x)-u(y))^2 \frac{c(x, y)}{|x-y|^{d+\alpha}} m(dy) \right) m(dx)=  f_u(x) m(dx)  .
$$
Then by (i), $\mu_{\<u\>}\in \K_\infty (X^{(1)})$ if
$f_u \in L^p(E; m)$ for some $p>d/\alpha$.
Note that $C^1_c(E)\subset \sF$.
 We next show that for $u\in C^1_c(E)$, $f_u\in L^p(E; m)$ for
 every $p\geq 1$.
  Clearly, by the mean value theorem,
 $$f_u (x)\leq \int_{\{y\in E: |y-x|< 1\}} \frac{c}{|x-y|^{d+\alpha-2}}
  m(dy) + \int_{\{y\in E: |y-x|\geq 1\}} \frac{ c }{|x-y|^{d+\alpha}} m(dy)
 $$
and so $f_u$ is bounded.
Let $K={\rm supp} [ u]$. Then for $x\in K^c$,
$$ f_u(x) =  \int_{K}   u(y)^2 \frac{c(x, y)}{|x-y|^{d+\alpha}} m(dy)
\leq \frac{c}{1+|x|^{d+\alpha}}.
$$
Since $(E, m)$ is a $d$-set, it follows that
$f_u\in L^p (E; m)$ for every $p\geq 1$. In particular,
we have $\mu_{\<u\>}\in \K_\infty (X^{(1)})$.

(iii) Suppose that $F$ is a bounded Borel function on
$E\times E$ satisfying \eqref{e:5.8a}.
Let $f(x)= \int_E |F(x, y)| N(x, y) m(dy)$. Then for every
$x\in E$,
\begin{eqnarray*}
 f(x)&\leq& \int_{\{y\in K: |y-x|< 1\}}
\frac{c}{|x-y|^{d+\alpha -\gamma}} m(dy) +
\int_{\{y\in K: |y-x|\geq 1\}}
\frac{c}{|x-y|^{d+\alpha}} m(dy)\\
& \leq & \frac{c}{d(x, K)^{d+\alpha}},
\end{eqnarray*}
where $d(x, K)$ denotes the Euclidean distance
between $x$ and $K$. It follows that
$f\in L^p(E; m)$ for every $p\geq 1$.
In particular, this implies that $F\in \J_\infty
(X^{(1)})$. \qed

\noindent{\bf Example 5.2.} {\bf (Symmetric diffusions)}
  Let $X$ be a symmetric diffusion
in  $\bR^n$, $n\geq 1$,
 with infinitesimal generator $\sL ={1\over 2} \sum_{i, j=1}^n
{\partial \over \partial x_i} \left( a_{ij} (x) {\partial \over
\partial x_j} \right)$, where matrix
$(a_{ij}(x))_{1\leq i, j\leq n}$ is uniformly elliptic and bounded,
that is, there is $\lambda >1$ such that for $m$-a.e. $x\in \bR^n$
and $\xi= (\xi_1, \cdots, \xi_n)\in \bR^n$,
$$ \lambda^{-1}\, \| \xi\|^2 \leq \sum_{i, j=1}^n
a_{ij}(x) \xi_i\xi_j \leq \lambda \, \| \xi \|^2.
$$
The Dirichlet form  $(\sE, \sF)$ in $L^2(\bR^n, dx)$ for $X$ is:
$\sF=W^{1, 2}(\bR^n)=\{ f\in L^2(\bR^n, dx): \ \nabla f\in L^2(\bR^n, dx)\}$
and $$ \sE (f, g)= {1\over 2} \int_{\bR^n} \sum_{i,j=1}^n
a_{ij}(x) {\partial f\over \partial x_i}
{\partial g \over \partial x_j} \, dx, \qquad f, g\in W^{1, 2}(\bR^n).
$$

\begin{lemma}\label{L:5.2} \begin{description}
\item{\rm (i)}
If $n\geq 3$, then  $ L^p(\bR^n; dx)\subset \K_\infty (X^{(1)})$
for every $p>n/2$.
When $n=1$ or $2$, then $L^p(\bR^n; dx)\subset \K_\infty (X^{(1)})$
for every $p\geq 1$.

\item{\rm (ii)} Suppose that $u\in L^2_{\rm loc}(\bR^n; dx)$
is bounded with
$\nabla u \in L^2(\bR^n; dx) \cap L^p (\bR^n; dx)$ for some $p>n$.
Then $u\in \sF_e$ with $\mu_{\<u\>}\in \K_\infty (X^{(1)})$.

\end{description}
\end{lemma}

\pf  (i) It is well known that the symmetric diffusion process $X$ has
a jointly H\"older continuous transition density function
$p(t, x, y)$ with respect to the Lebesgue measure $m(dx):=dx$ on $\bR^n$.
Moreover, $p(t, x, y)$ enjoys the following celebrated Anroson's estimate:
there are constants $c_1, c_2\geq 1$ so that
for every $t>0$ and $x, y\in \bR^n$,
\begin{equation}\label{e:5.15}
c_1^{-1} t^{-n/2} e^{-c_2|x-y|^2/t} \leq p(t, x, y)
\leq c_1 t^{-n/2} e^{- |x-y|^2/ (c_2t)} .
\end{equation}
Note that
\begin{eqnarray*}
\int_0^\infty e^{-t} t^{-n/2} e^{- r^2/(c_2)t} dt
&  \stackrel{t=r^2 u}{=} &
 r^{2-n} \int_0^\infty u^{-n/2} e^{-1/(c_2u)} e^{-r^2 u} du\\
 &\leq & c_3 r^{2-n} \int_0^\infty (u^{-n/2-1}\wedge 1) e^{-1/(2c_2u)-r^2u/2}du \\
 &\leq & c_3 r^{2-n} e^{-c_4r}\int_0^\infty (u^{-n/2-1}\wedge 1)  du
  =c_5 r^{2-n} e^{-c_4r}.
\end{eqnarray*}
Thus we have by \eqref{e:5.15} and the above that
\begin{equation}\label{e:5.16}
G_1(x, y)=\int_0^\infty e^{-t}p(t, x, y)dt
\leq c_1c_5 |x-y|^{2-n}  e^{ -c_4 |x-y|}.
\end{equation}
Just as in the proof of Lemma \ref{L:5.1},
using  H\"older inequality and \eqref{e:5.16}, one can show that when $n\geq 3$,
$L^p(\bR^n;dx)\subset \K_\infty (X^{(1)})$ for
every $p>n/2$. When $n=1$ or  $2$, $ L^p(\bR^n;dx)\subset
 \K_\infty (X^{(1)})$ for every $p\geq 1$.

(ii) It is known (see, e.g., Example 1.5.2 of \cite{FOT})
that the extended Dirichlet space
$$
\sF_e= \left\{ f\in L^2_{\loc}(\bR^n; dx): \ \nabla f\in L^2(\bR^n; dx)
 \right\}.
$$
By \eqref{e:3.1a}, for bounded $u\in \sF_e$, its energy measure is $\mu_{\<u\>}(dx)
= \sum_{i, j=1}^n a_{ij}(x) {\partial u \over \partial x_i}
{\partial u \over \partial x_j}  dx$.
Thus by (i) above, a bounded locally $L^2$-integrable function
$u$ with $\nabla u \in L^2(\bR^n; dx) \cap L^p(\bR^n; dx)$
for some $p>n$ is a function in $\sF_e$ with $\mu_{\<u\>}$
in the Kato class of $X$. \qed

\bigskip

We refer the reader to \cite[Examples 2.2 and 2.3]{C2}
for examples of Kato classes $\K_\infty (X)$ and
$\J_\infty (X)$ when $X$ is a symmetric $\alpha$-stable process,
respectively, a censored $\alpha$-stable process, in a bounded
$C^{1,1}$-open set.

\vskip 0.3truein

\begin{singlespace}

\end{singlespace}

\vskip 0.3truein

{\bf Zhen-Qing Chen}

Department of Mathematics, University of Washington, Seattle,
WA 98195, USA

E-mail: \texttt{zqchen@uw.edu}

\end{document}